\documentclass[reqno,centertags, 12pt]{amsart}
\usepackage{amsmath,amsthm,amscd,amssymb}
\usepackage{latexsym}   

\newcommand{\bbR}{{\mathbb{R}}}

\newcommand{\bbC}{{\mathbb{C}}}

\newcommand{\calH}{{\mathcal H}}

\newcommand{\dott}{\,\cdot\,}

\newcommand{\lb}{\label}
\newcommand{\f}{\frac}

\newcommand{\loc}{\text{\rm{loc}}}

\newcommand{\ac}{\text{\rm{ac}}}
\newcommand{\singc}{\text{\rm{sc}}}
\newcommand{\pp}{\text{\rm{pp}}}
\newcommand{\sing}{\text{\rm{sing}}}

\newcommand{\bi}{\bibitem}

\newcommand{\beq}{\begin{equation}}
\newcommand{\eeq}{\end{equation}}
\newcommand{\ba}{\begin{align}}
\newcommand{\ea}{\end{align}}
\newcommand{\veps}{\varepsilon}
\DeclareMathOperator{\Real}{Re}

\allowdisplaybreaks
\numberwithin{equation}{section}
\newtheorem{theorem}{Theorem}[section]
\newtheorem{proposition}[theorem]{Proposition}
\newtheorem{lemma}[theorem]{Lemma}
\newtheorem{corollary}[theorem]{Corollary}
\theoremstyle{definition}

\theoremstyle{remark}

\newcommand{\abs}[1]{\lvert#1\rvert}

\begin{document}
\title[Stability of Singular Spectral Types]{Stability of Singular \\
Spectral Types under Decaying Perturbations}
\author[A. Kiselev, Y. Last, and B. Simon]{Alexander Kiselev$^1$, Yoram Last$^2$, 
and Barry Simon$^3$}
\date{October 8, 2001}
\dedicatory{Dedicated to Jean Michel Combes on his 60th birthday}

\footnotetext[1]{ Supported in part by NSF Grant 
No.~DMS-0102554.}
\footnotetext[2]{ Supported in part by The Israel Science Foundation 
Grant No.~447/99 and by an Allon fellowship.}
\footnotetext[3]{ Supported in part by NSF Grant No.~DMS-9707661.}

\address{A. Kiselev, Department of Mathematics, University of Chicago, Chicago, 
IL 60637, USA}
\email{kiselev@math.uchicago.edu}

\address{Y. Last, Institute of Mathematics, The Hebrew University, 91904 
Jerusalem, Israel}
\email{ylast@math.huji.ac.il}

\address{B. Simon,  Division of Physics, Mathematics, and Astronomy, 
253-37, California Institute of Technology, Pasadena, CA~91125, USA} 
\email{bsimon@caltech.edu}

\begin{abstract} We look at invariance of a.e.~boundary condition spectral behavior 
under perturbations, $W$, of half-line, continuum or discrete Schr\"odinger operators. 
We extend the results of del Rio, Simon, Stolz from compactly supported $W$'s to 
suitable short-range $W$. We also discuss invariance of the local Hausdorff dimension 
of spectral measures under such perturbations.
\end{abstract}
\maketitle

\section{Introduction} \lb{s1}

\bigskip
We want to discuss aspects of the spectral theory of Schr\"odinger operators on a half-line, 
both continuous
\begin{equation} \lb{1.1}
(Hu)(x) = -\f{d^2}{dx^2} + V(x)
\end{equation}
on $L^2 (0,\infty; dx)$ and discrete
\begin{equation} \lb{1.2}
(hu)(n) = u(n+1) + u(n-1) + V(n) u(n)
\end{equation}
on $\ell^2 (\{1,2,\dots\})$ with $u(0)$ determined by the boundary condition. These 
operators have a boundary condition determined by a parameter $\theta$ in $[0,\pi)$:
\begin{equation} \lb{1.3}
u(0) \cos(\theta) + u'(0) \sin(\theta) =0
\end{equation}
in the continuum case and
\begin{equation} \lb{1.4}
u(0) \cos \theta + u(1) \sin(\theta) = 0
\end{equation}
in the discrete case. Thus \eqref{1.4} is equivalent to defining
\[
(h_\theta u)(1) = u(2) + [V(1) - \tan(\theta)] u(1).
\]

In some places below, we will suppose $V(x)$ is bounded in the continuum case 
for reasons that will become clear. In the discrete case, we will need 
boundedness only once.

We will use $H_\theta$ and $h_\theta$ to indicate the operators with boundary 
condition. It is well known (see, e.g., Simon \cite{van}) that there are 
spectral measures $d\rho_\theta (\lambda)$ for $H_\theta$ and $h_\theta$ (so 
that $H_\theta$ or $h_\theta$ is unitarily equivalent to multiplication by 
$\lambda$ on $L^2 (\bbR, d\rho_\theta (\lambda))$) normalized so that 
\begin{equation} \lb{1.5}
\int_{\theta=0}^\pi d\rho_\theta (\lambda) \f{d\theta}{\pi} = d\lambda.
\end{equation}

A major theme in this paper (as in many recent papers) is the relation of spectral 
properties with solutions of the differential/difference equation. Given $V$ and 
$\theta$, for each $\lambda\in\bbC$, we will define $\varphi_{1,\theta}(\lambda, x)$ 
(or $\varphi_{1,\theta} (\lambda,n)$) to be the solution of
\begin{equation} \lb{1.6}
H\varphi = \lambda\varphi \qquad (\text{or } h\varphi=\lambda\varphi)
\end{equation}
(intended as a differential/difference equation with no $L^2$ condition at 
$\infty$) obeying the boundary condition \eqref{1.3}/\eqref{1.4} and normalized 
by
\begin{equation} \lb{1.7}
\varphi_{1,\theta}(\lambda, 0) = \sin(\theta) \qquad \varphi'_{1,\theta} 
(\lambda, 0) = -\cos(\theta)
\end{equation}
(or $\varphi_{1,\theta}(\lambda, 1) = -\cos(\theta)$ in the discrete case). We 
will also define
\[
\varphi_{2,\theta} \equiv \varphi_{1, (\theta -\pi/2)}.
\]

While we consider $\theta\in [0,\pi)$ in the basic definition of $\varphi_{1,\theta}$, 
it makes sense for all $\theta$ with $\varphi_{1,\theta + n\pi} = (-1)^n \varphi_{1,\theta}$. 
In particular, in the last equation $\theta - \pi/2$ lies in $[-\pi/2, \pi/2)$. With this 
definition, the Wronskian obeys
\begin{equation} \lb{1.8}
W(\varphi_{1,\theta}, \varphi_{2,\theta}) =1
\end{equation}
with $W(f,g) = fg' - f'g$ in the continuum case and $W(f,g)(n) = f(n) g(n+1) - 
f(n+1) g(n)$ in the discrete case.

Following Jitomirskaya-Last \cite{JL}, for $L>0$, we define
\[
\|f\|^2_L = \int_0^L |f(x)|^2\, dx 
\]
in the continuum case and
\[
\|f\|^2_L = \sum_{n=1}^{[L]} |f(n)|^2 + (L - [L]) |f([L] + 1)|^2
\]
in the discrete case (so $\|f\|^2_L$ is the obvious analog at integer $L$, with 
linear interpolation in between).

When one looks at the decomposition of $d\rho_\theta$ into spectral types, for example, 
into a.c., s.c., and pure point pieces (see Reed-Simon \cite{RS1}), a basic pair of 
facts says that the a.c.~spectrum is stable and the singular spectrum is unstable --- 
explicitly (see Simon \cite{van} for references), the essential support of 
$d\rho^{\ac}_\theta$ is $\theta$ independent, while for any pair $\theta \neq \theta'$, 
$d\rho^{\sing}_\theta$ and $d\rho^{\sing}_{\theta'}$ are mutually singular. These 
facts seem to be at variance with the notion that spectral properties should  
depend on the behavior of $V$ at infinity since they suggest that $d\rho^{\sing}$ 
will be unstable under perturbations of compact support. The resolution of this 
conundrum is the idea of del Rio, Simon, and Stolz \cite{DRSS} that one should 
look at the union over $\theta$ of spectral supports. Explicitly, we proceed as 
follows:

\smallskip
\noindent{\bf Definition (Gilbert-Pearson \cite{GP}).} We say there is a subordinate 
solution at energy $\lambda\in\bbR$ if and only if there is some $\theta\in [0,\pi)$ 
so $\lim_{L\to\infty} \|\varphi_{1,\theta}\|_L / \|\varphi_{2,\theta}\|_L = 0$. 
$\theta$ is necessarily unique and we call it $\theta(\lambda)$.

\smallskip
\noindent{\bf Definition.} 
\begin{align*}
P &= \{\lambda \mid \varphi_{1,\theta(\lambda)}\in L^2\} \\
S &= \{\lambda\mid\text{there is a subordinate solution but }\varphi_{1,\theta(\lambda)}
\notin L^2\} \\
L &= \{ \lambda \mid\text{ there is no subordinate solution}\}
\end{align*}
When we need to discuss the $V$-dependence of these sets, we will write $P(V)$, etc.

Then:

\begin{theorem} \lb{T1.1} 
\begin{enumerate}
\item[(i)] $P = \cup_\theta \, \sigma_{\pp} (H_\theta)$. 
\item[(ii)] $L = \text{essential support of $\sigma_{\ac} (H_\theta)$ 
for all $\theta$}$.
\item[(iii)] For any $\theta$, $d\rho^{\singc}_\theta = d\rho_\theta (S \cap \cdot)$ 
and if $\tilde S$ is any other set with that property, then $|S \triangle \tilde S| 
=0$ where $|\dott|$ is Lebesgue measure.
\end{enumerate}
\end{theorem}

\noindent{\it Remarks.} 1. This is close to a theorem in \cite{DRSS}, although 
$S$ and $L$ are defined differently there.

\smallskip
2. $\sigma_{\pp}$ in (i) means the set of eigenvalues, not their closure.

\smallskip
3. (i) is obvious since $\lambda\in\sigma_{\pp} (H_{\theta(\lambda)})$ if and only if 
$\varphi_{1,\theta(\lambda)}(\lambda, \cdot)\in L^2$.

\smallskip
4. (ii) is the main result of Gilbert-Pearson \cite{GP}.

\smallskip
5. That $d\rho^{\singc}_\theta(P) =0$ is obvious since $d\rho^{\sing}_\theta$ is 
mutually singular to each $d\rho^{\pp}_{\theta'}$ for $\theta'\neq \theta'$ and 
$d\rho^{\singc}_\theta$ is obviously mutually singular to $d\rho^{\pp}_\theta$.

\smallskip
6. That $d\rho^{\singc}_\theta (L)=0$ is a result of Gilbert-Pearson showing that 
$d\rho^{\singc}_\theta = d\rho_\theta (S\cap\cdot)$.

\smallskip
7. The $|S\triangle \tilde S|=0$ result follows from \eqref{1.5}.

\medskip
Since $P$, $L$, $S$ are defined purely in terms of the behavior of solutions at 
infinity, the following result of del Rio et al.~\cite{DRSS} is immediate:

\begin{theorem} \lb{T1.2} Let $V=V_0 + W$ where $W$ has compact support. Then
$P(V) = P(V_0)$, $L(V)=L(V_0)$, $S(V)=S(V_0)$.
\end{theorem}

A major theme of this paper will be to examine when this result still holds for 
$W$'s not of compact support. Before discussing our theorems, we will further refine 
the set $S$ in connection with the breakdown of singular spectrum according to 
Hausdorff measures and dimensions.

As usual for $\alpha\in (0,1)$, $\alpha$-dimensional Hausdorff measure is defined on 
Borel sets, $T$, by
\[
h^\alpha (T) \equiv \lim_{\delta\to 0} \inf_{\delta-\text{covers}}\sum_{\nu=1}^\infty 
|b_\nu|^\alpha, 
\]
where a $\delta$-cover is a countable collection of intervals each of length at most 
$\delta$ so $T\subset \cup_{\nu=1}^\infty b_\nu$. $h^1$ is Lebesgue measure and $h^0$ 
is counting measure.

Given $\alpha\in [0,1]$ (following Rogers and Taylor \cite{RT1,RT2}; see also Last 
\cite{La}), we define a measure $\mu$ to be $\alpha$-continuous ($\alpha$c) if 
$\mu(S)=0$ for any set $S$ with $h^\alpha (S) =0$ and $\alpha$-singular ($\alpha$s) 
if it is supported on a set of $S$ with $h^\alpha (S) =0$. For every such $\alpha$ 
and any measure $\mu$, one can uniquely decompose $\mu = \mu^{\alpha\text{c}} + 
\mu^{\alpha\text{s}}$ with $\mu^{\alpha\text{c}}$ $\alpha$-continuous and 
$\mu^{\alpha\text{s}}$, $\alpha$-singular.

We call a measure zero-dimensional if it is supported on a set $S$ with $h^\alpha (S) =0$ 
for all $\alpha >0$. We call it one-dimensional if it is 
$\alpha$-continuous for all $\alpha <1$.

It will be useful, following Jitomirskaya-Last, to have a pair of inverse functions 
$A,B : [0,1]$ to $[0,1]$ by
\begin{align*}
B(\alpha) &= \alpha /(2-\alpha) \\
A(\beta) &= 2\beta / (1+\beta).
\end{align*}

\smallskip
\noindent{\bf Definition.} Let $\lambda\in S$, the set of energies for which there is a 
non-$L^2$ subordinate solution. Define
\[
\beta(\lambda) = \liminf_{L\to\infty} \big[ \ln \|\varphi_{1,\theta(\lambda)}\|_L \big/ 
\ln \|\varphi_{2,\theta(\lambda)}\|_L \big].
\]

Notice that since $\varphi_{1,\theta}\notin L^2$, $\|\varphi_{1,\theta(\lambda)}\|_L 
\to\infty$ as $L\to\infty$ and since $\varphi_1$ is subordinate, eventually 
$\|\varphi_{2,\theta}\|_L \geq \|\varphi_{1,\theta}\|_L$, and thus
\begin{equation} \lb{1.9}
\lim_{L\to\infty} \| \varphi_{2,\theta}\|_L = \infty
\end{equation}
and
\[
0\leq\beta (\lambda) \leq 1.
\]

When we want to indicate the $V$-dependence of $\beta$, we will write $\beta 
(\lambda; V)$. We note the following elementary:

\begin{proposition} \lb{P1.3} If $\beta > \beta(\lambda)$, then
\begin{equation} \lb{1.10}
\varliminf \| \varphi_1\|_L \big/ \|\varphi_2 \|^\beta_L = 0
\end{equation}
and if $\beta < \beta (\lambda)$,
\begin{equation} \lb{1.11}
\lim \|\varphi_1 \|_L \big/ \|\varphi_2 \|^\beta_L = \infty.
\end{equation}
\end{proposition}

\begin{proof} Write
\[
\|\varphi_1\|_L \big/ \|\varphi_2\|^\beta_L = \exp \bigg[ \ln \|\varphi_2\|_L 
\bigg\{ \f{\ln \|\varphi_1\|_L}{\ln \|\varphi_2\|_L} - \beta \bigg\}\bigg].
\]
By \eqref{1.9}, $\ln \|\varphi_2\|_L\to\infty$. If $\beta >\beta (\lambda)$, 
then there is a subsequence where the expression in $\{ \,\, \}$ goes to 
$\beta (\lambda) -\beta <0$, so a subsequence where the expression in 
$[\,\, ]$ goes to $-\infty$ and \eqref{1.10} holds. If $\beta <\beta 
(\lambda)$, then eventually the expression in $\{\,\,\}$ is larger than 
$\f12 (\beta(\lambda)-\beta)$, and so \eqref{1.11} holds.
\end{proof}

For each $\beta_0$, decompose $S$ into four sets:
\begin{align*}
S^{++}_{\beta_0} &= \{ \lambda \mid \beta_0 > \beta (\lambda)\} \\
S^{--}_{\beta_0} &= \{ \lambda \mid \beta_0 < \beta (\lambda) \}\\
S^+_{\beta_0} &= \{\lambda \mid\beta_0 = \beta (\lambda) \text{ and \eqref{1.10} 
holds for } \beta_0 = \beta (\lambda) \} \\
S^-_{\beta_0} &= \{\lambda \mid \beta_0 = \beta (\lambda) \text{ and } \varliminf 
\|\varphi_1\|_L / \|\varphi_2\|^{\beta_0} > 0\}.
\end{align*}
Thus \eqref{1.10} holds for $\beta =\beta_0$ if and only if $\lambda\in S^{++}_{\beta_0} 
\cup S^+_{\beta_0}$.

It follows from Theorem~\ref{T1.1} and the discussion following equation~(2.2) of 
Jitomirskaya-Last \cite{JL} that

\begin{theorem} \lb{T1.4} Let $\beta_0 = B (\alpha_0)$.
\begin{enumerate}
\item[(i)] $d\rho^{\alpha_0\text{\rm{c}}}_\theta = d\rho_\theta ((S^-_{\beta_0} 
\cup S^{--}_{\beta_0} \cup L) \cap \cdot )$
\item[(ii)] $d\rho^{\alpha_0\text\rm{{s}}}_\theta = d\rho_\theta ( (S^+_{\beta_0} 
\cup S^{++}_{\beta_0} \cup P) \cap \cdot)$
\item[(iii)] $d\rho^{\singc}_\theta$ is one-dimensional for a.e.~$\theta$ if and only 
if $\beta =1$ a.e.~on $S$.
\item[(iv)] $d\rho^{\singc}_\theta$ is zero-dimensional for a.e.~$\theta$ if and only 
if $\beta =0$ a.e.~on $S$.
\end{enumerate}
\end{theorem}

\noindent{\it Remark.} More generally, $d\rho^{\singc}_\theta$ has exact dimension 
$\alpha_0$ for a.e.~$\theta$ if $\beta = B (\alpha_0)$ for a.e.~$\lambda\in S$.

\medskip

Clearly, $\beta$ only depends on $V$ near infinity, so we extend the result of 
del Rio et al. \cite{DRSS} to handle dimensional decomposition of $d\rho$ via

\begin{theorem} \lb{T1.5} Let $V=V_0 +W$ where $W$ has compact support. Then 
$\beta (\lambda;V) = \beta (\lambda; V_0)$.
\end{theorem}

The purpose of this paper is to study when invariance results of the genre of 
Theorems~\ref{T1.2} and \ref{T1.5} extend to cases where $W$ does not have 
compact support but has ``suitable" decay; that is, we want to determine what 
suitable decay is. For the a.c.~spectrum, the standard rate of decay is 
$W\in L^1$:

\begin{theorem} \lb{T1.6} In the continuum case, suppose $V_0$ and $V\equiv 
V_0 + W$ are such that $H_0 + V_0$ and $H_0 +V$ are bounded below by $\veps 
H_0 -c$. In the discrete case, no hypothesis is needed on $V_0$. Suppose 
that $W\in L^1$ {\rm(}or $\ell^1{\rm)}$. Then
\begin{equation} \lb{1.12}
\abs{L(V)\triangle L(V_0)}=0.
\end{equation}
\end{theorem}

\begin{proof} In the discrete case, $W$ is trace class, and in the continuum case, 
$(H_0 +1)^{-1/2} W(H_0 +1)^{-1/2}$ is trace class. So $(H_0 + V + c-1)^{-1} -
(H_0 + V_0 + c + 1)^{-1}$ is trace class. The trace class theory of scattering 
\cite{RS3} implies that $H_0$ on $\calH_{\ac} (H_0)$ is unitarily equivalent to $H_0 
+ V_0$ on $\calH_{\ac}(H_0 + V_0)$ from which \eqref{1.12} follows by 
Theorem~\ref{T1.1}.
\end{proof}

\noindent{\it Remark.} We conjecture that \eqref{1.12} holds if $W$ is merely 
assumed in $L^2$. In \cite{KLS}, we made this conjecture when $V_0 =0$ and it was 
proven by Deift-Killip \cite{DK}. Killip \cite{Ki} proved the result when 
$V_0$ is periodic. We conjecture the result for all $V_0$.

\medskip

We now turn to the substantially new results in this paper. As spectrum moves from 
the most smooth (a.c.) to the least smooth (point), we need to successively 
strengthen the conditions on the perturbation $W$.

We begin with several results we prove in Section~\ref{s3} concerning point 
spectrum that all hold in the discrete and continuum case.

\begin{theorem} \lb{T1.7} For each $\lambda\in P(V_0)$, define
\begin{equation} \lb{1.13}
f_+ (\lambda, x)=(1+ \abs{x}) \sup_{\abs{y}\leq x} \, 
\abs{\varphi_{2,\theta (\lambda)} (y)}.
\end{equation}
Suppose that for all $\lambda\in Q \subseteq P(V_0)$, we have that
\[
\int \abs{W(x)} f_+ (\lambda, x)\, dx <\infty 
\]
and that the $L^2$ solution is bounded. Then $Q\subseteq P(V_0 +W)$.
\end{theorem}

\noindent{\it Remarks.} 1. In \eqref{1.13}, one can replace $(1 +\abs{x})$ by 
$(1 + \abs{x})^\gamma$ for any $\gamma >\f12$.

2. By a Sobolev estimate if $\varphi, \varphi'\in L^2$, then $\varphi\in L^\infty$, 
so, for example, if $V_0$ is bounded from below, $L^2$ solutions will be bounded.

\medskip
When $V_0$ is bounded, $f_+$ does not grow faster than exponentially for any 
$\lambda$.

\begin{corollary} \lb{C1.8} Let $V_0$ be bounded and suppose that
\[
\int \abs{W(x)}\, e^{A|x|}\, dx <\infty
\]
for all $A>0$. Then
\[
P(V_0) = P(V_0 +W).
\]
\end{corollary}

Finally, we have a result on preservation of Lyapunov behavior. Recall that we say 
there is Lyapunov behavior at energy $\lambda$ if the transfer matrix
\[
T_\lambda (0,x) = 
\begin{pmatrix}
\varphi'_{1,\theta}(x) & \varphi'_{2,\theta}(x) \\
\varphi_{1,\theta}(x) & \varphi_{2,\theta}(x) 
\end{pmatrix}
\]
obeys
\begin{equation} \lb{1.14}
\lim_{x\to\infty} \f{1}{\abs{x}}\, \ln \|T_\lambda (0,x)\| \equiv \gamma (\lambda).
\end{equation}

\begin{theorem} \lb{T1.9} Suppose $V_0$ has Lyapunov behavior at energy $\lambda$ 
and that for some $\veps >0$,
\[
\int \abs{W(x)} \, e^{\veps |x|}\, dx <\infty.
\]
Then $V_0 +W$ has Lyapunov behavior at $\lambda$ with the same value of $\gamma$.
\end{theorem}

\noindent{\it Remarks.} 1. If $\gamma >0$, we have much more than merely the same 
Lyapunov behavior.

2. Theorem~\ref{T1.9} isn't new. It is essentially a special case of Theorem
4.I of \cite{Ru}.

\medskip
In Section~\ref{s4}, we will discuss stability of singular spectrum and its components. 
Our results will hold only for energies with an extra condition.

\smallskip
\noindent{\bf Definition.} An energy $\lambda$ is called regular if and only if for some 
$\theta$ ($=\theta (\lambda)$ if there is a subordinate solution) we have for all $\veps 
>0$,
\begin{equation} \lb{1.15}
\| \varphi_{1,\theta}\|_L \leq C_\veps L^{1/2 +\veps}.
\end{equation}

\smallskip
By the general theory of eigenfunction expansions \cite{Ber,ssg}, a.e.~$\lambda$ is regular both 
with respect to each $d\rho_\theta$, and so by \eqref{1.5} for a.e.~$\lambda$ with 
respect to Lebesgue measure $d\lambda$. Indeed, we could replace $L^{1/2 +\veps}$ by 
$L^{1/2}(\ln L)^\kappa$ for any $\kappa >\f12$.

\smallskip
\noindent{\it Remark.} If $V(x)= -\f{3}{16}\, x^{-2}$ for large $x$, then the subordinate 
solution at $\lambda =0$ is $\sim x^{1/4}$ at infinity. So $\|\varphi_1\|_L \sim L^{3/4}$ 
and $\lambda =0$ is not a regular energy, so not all energies need to be regular.

\medskip 
In the discrete case, constancy of the Wronskian implies
\begin{equation} \lb{1.16}
\|\varphi_{1,\theta}\|_L \, \|\varphi_{2,\theta}\|_L \geq \tfrac12 (L-1),
\end{equation}
but in the continuum case, this is not automatic since the Wronskian involves $\varphi'$. 
But, by a Sobolev estimate, if $V$ is bounded (uniform locally $L^1$ will do!), then
\begin{equation} \lb{1.17}
\|\varphi_{1,\theta}\|_L \, \|\varphi_{2,\theta}\|_L \geq c(L-1)
\end{equation}
for some $c$, dependent on $V$ and $\lambda$, and so we will need to suppose that $V$ 
is bounded in the continuum case.

\smallskip
\noindent{\it Remark.} The case $V(x) =-x$ where $\|\varphi_{1,\theta}\|_L \sim 
\|\varphi_{2,\theta}\|_L \sim L^{1/4}$ shows \eqref{1.17} really can fail if $V$ is 
unbounded.

\medskip
Here are the theorems we will prove in Section~\ref{s4}.

\begin{theorem} \lb{T1.10} In the continuum case, suppose $V_0$ is bounded. Let 
$\lambda \in S(V_0)$ be a regular energy with $\beta (\lambda, V_0)=1$. Suppose that
\begin{equation} \lb{1.18}
\abs{W(x)} \leq C(1+\abs{x})^{-1-\veps}
\end{equation}
for some $\veps >0$. Then $\lambda\in S(V_0 +W)$ with $\beta (\lambda, V_0 +W) =1$. In 
particular, if, for $V_0$, $H_\theta$ has one-dimensional spectrum for a.e.~$\theta$, 
the same is true for $V_0 +W$.
\end{theorem}

\begin{theorem} \lb{T1.11} In the continuum case, suppose $V_0$ is bounded. Let 
$\lambda\in S(V_0)$ be a regular energy. Suppose that for all $\eta >0$, 
\begin{equation} \lb{1.19}
\abs{W(x)} \leq C_\eta (1+\abs{x})^{-\eta}.
\end{equation}
Suppose that $\beta (\lambda, V_0)\neq 0$. Then $\lambda\in S(V_0 +W)$ and 
$\beta (\lambda, V_0 +W) =\beta (\lambda, V_0)$. Suppose $\beta (\lambda, V_0) =0$. 
Then either $\lambda\in S(V_0 +W)$ with $\beta (\lambda, V_0 +W) =0$ or $\lambda 
\in P(V_0 +W)$.
\end{theorem}

\noindent{\it Remarks.} 1. The latter shows that having zero-dimensional spectrum is 
preserved under perturbations obeying \eqref{1.19}, although to preserve point 
spectrum, we need a stronger exponential bound.

\smallskip
2. In fact, our proof shows that for a given $\beta (\lambda, V_0)=\beta_0$, we only 
need \eqref{1.19} for some
\[
\eta > \f{1}{\beta_0}\, .
\]
In terms of the case of Hausdorff dimension $\alpha$, one needs
\begin{equation} \lb{1.20}
\eta > \f{2}{\alpha } -1.
\end{equation}

\medskip
We will prove our new results, Theorem~\ref{T1.7}, Corollary~\ref{C1.8}, and 
Theorems~\ref{T1.9}--\ref{T1.11}, by proving stability of the asymptotics of 
solutions of the Schr\"odinger differential/difference equation. We use $\varphi_-$ 
for $\varphi_{1,\theta (\lambda)}$, the subordinate solution with potential $V_0$,  
and $\varphi_+$ for $\varphi_{2,\theta(\lambda)}$. The basic construction we will 
use is variation of parameters. That is, we will write (in the continuum case):
\begin{align} 
\psi(x) &= u_1 (x) \varphi_-(x) + u_2(x) \varphi_+(x) \lb{1.21} \\
\psi' (x) &= u_1(x) \varphi'_-(x) + u_2 (x) \varphi'_+(x). \lb{1.22}
\end{align}
With $u(x) = \binom{u_1 (x)}{u_2(x)}$, the differential equation for $\psi$ 
is equivalent, given the normalization \eqref{1.7}, to
\begin{equation} \lb{1.23}
u'(x) = A(x) u(x)
\end{equation}
with 
\begin{equation} \lb{1.24}
A(x) = -W(x) 
\begin{pmatrix} 
\varphi_+ (x) \varphi_-(x) & \varphi_+(x)^2 \\
-\varphi_- (x)^2 & -\varphi_+(x) \varphi_-(x)
\end{pmatrix} .
\end{equation}
\eqref{1.24} is sometimes written (e.g., in \cite{JL}) in the integral form:
\begin{equation} \lb{1.25}
\begin{split}
\psi(x) &= u_1(x_0) \varphi_-(x) + u_2 (x_0) \varphi_+ (x) \\
&\qquad - \int_{x_0}^x 
W(y) [\varphi_+(x) \varphi_-(y) - \varphi_-(x) \varphi_+(y)] \psi (y)\, dy.
\end{split}
\end{equation}

In the discrete case, the result is similar. One writes
\begin{align}
\psi(n) &= u_1(n) \varphi_-(n) + u_2 (n) \varphi_+(n) \lb{1.26} \\
\psi (n-1) &= u_1 (n) \varphi_-(n-1) + u_2(n) \varphi_+ (n-1). \lb{1.27}
\end{align}
\eqref{1.23} becomes
\begin{equation} \lb{1.28}
u(n+1) - u(n) = A(n) u(n), 
\end{equation}
where
\begin{equation} \lb{1.29}
A(n) = -W(n) \begin{pmatrix}
\varphi_+(n) \varphi_-(n) & \varphi_+ (n)^2 \\
-\varphi_- (n)^2 & -\varphi_+ (n) \varphi_+ (n) \varphi_-(n)
\end{pmatrix}
\end{equation}
or its integral form
\begin{equation} \lb{1.30}
\begin{split}
\psi(n) &= u_1 (n_0) \varphi_-(n) + u_2 (n_0) \varphi_+(n) \\
&\qquad + \sum_{j=n_0 +1}^n W(j) [\varphi_+ (n) \varphi_- (j) - \varphi_-(n) 
\varphi_+(j)] \psi(j).
\end{split}
\end{equation}

The standard control for perturbing solutions at infinity is to require 
$\int_{x_0}^\infty \|A(x)\|\, dx <\infty$. For the diagonal matrix elements 
of $A$, that cannot be improved without detailed oscillation estimates, but 
it is well known that one can try to trade off the growth of one off-diagonal 
matrix element by the decay of the other. In Section~\ref{s2}, we present a 
version of this fact made for our applications. These ideas are not new; 
for example, our method of proof is patterned after problem XI.97 of Reed-Simon 
\cite{RS3}. In Section~\ref{s3}, we present the results of stability of a 
solution $L^2$ at $\infty$ and in Section~\ref{s4}, the results on stability of 
polynomially bounded solutions. The Appendix discusses some results concerning
the preservation of WKB asymptotic behavior of solutions.

%%%%%%%%%%%%%%%%%%%%%%%%%%%%%%%%%%%%%%%%%%%%%%%%%%%%%%%%%%
\section{A Perturbation Lemma} \lb{s2}
%%%%%%%%%%%%%%%%%%%%%%%%%%%%%%%%%%%%%%%%%%%%%%%%%%%%%%%%%%

In this preliminary section, we will be interested in solutions of
\begin{equation} \lb{2.1}
u'(x) = A(x) u(x),
\end{equation}
where
\begin{equation} \lb{2.2}
A(x) = \begin{pmatrix} 
a_{11}(x) & a_{12}(x) \\
a_{21}(x) & a_{22}(x) 
\end{pmatrix}
\end{equation}
is in $L^1_{\loc} [0,\infty)$ and
\begin{equation} \lb{2.3}
u(x)= \binom{u_1(x)}{u_2(x)}
\end{equation}
is a two-component vector. By a solution of \eqref{2.1}, we mean an absolutely continuous 
function so that \eqref{2.1} holds for a.e.~$x$. As usual, given any $x_0$ and 
$\omega\in \bbC^2$, there is a unique solution of \eqref{2.1} with $u(x_0) = \omega$. 

We will use a pair of non-negative functions $f_\pm (x)$ with
\begin{equation} \lb{2.4}
f_+ (x) f_-(x) \geq 1
\end{equation}
and $f_+$ monotone increasing and $f_-$ monotone decreasing (in some applications, we 
will take $f_\pm = e^{2(\pm \gamma + \veps)|x|}$ so you can have this example in 
mind). Define
\begin{equation} \lb{2.5}
G(x) = \max (|a_{11}(x)| + |a_{12}(x)| f_-(x), |a_{21}(x)| f_+(x) + 
|a_{22}(x)|).
\end{equation}

\begin{lemma} \lb{L2.1} Define $\| \dott\|^\pm_x$ as norms on $\bbC^2$ by
\begin{align*}
\| \omega \|^+_x &= \max (|\omega_1|, |\omega_2| f_+ (x)) \\
\| \omega \|^-_x &= \max (f_- (x) |\omega_1|, |\omega_2|).
\end{align*}
Then
\begin{equation} \lb{2.6}
\| A(x) \omega\|^\pm_x \leq G(x) \|\omega\|^\pm_x.
\end{equation}
\end{lemma}

\begin{proof} We will prove the $\| \dott\|^+$ result. The $\| \dott\|^-$ is similar. 
Note that 
\begin{align*}
|(A(x)\omega)_1| &\leq |a_{11}(x)| \, |\omega_1| + |a_{12}| f_+ (x)^{-1} f_+ 
(x)|\omega_2| \\
&\leq [|a_{11}(x)| + |a_{12}| f^{-1}_+(x) ] \, \|\omega\|^+_x \leq G(x) \|\omega\|^+_x
\end{align*}
since $f^{-1}_+ \leq f_-$ by \eqref{2.4} and
\begin{align*}
f_+ (x) |(A(x)\omega)_2| &\leq [|a_{21}(x) f_+(x) |\omega_1| + |a_{22}(x)| f_+(x) |\omega_2|] \\
&\leq (|a_{21}(x)| f_+(x) + |a_{22}(x)|) \|\omega\|^+_x \\
&\leq G(x) \|\omega\|^+_x.
\end{align*} 
\end{proof}

\begin{theorem} \lb{T2.2} Suppose $f_+$ is monotone increasing, $f_-$ is monotone decreasing, 
\eqref{2.4} holds, and 
\[
\int_x^\infty G(y)\, dy < \infty.
\]
Then there exist solutions $u^\pm$ of \eqref{2.1} so that as $x\to\infty$
\begin{enumerate}
\item[(i)] $u^-_1(x)\to 1, \quad f_+ (x) u^-_2(x) \to 0$
\item[(ii)] $u^+_1 (x) f_- (x) \to 0, \quad  u^+_2 (x) \to 1$.
\end{enumerate}
\end{theorem}

\begin{proof} Define $u^{-(n)}$ by
\begin{align*}
u^{-(0)} &= \binom{1}{0} \\
u^{-(n+1)}(x) &= -\int_x^\infty A(y) u^{-(n)} (y)\, dy,
\end{align*}
where we will deal with the convergence of the integral below. Since $f_+$ is increasing,  
if $y > x$, then $\|\omega\|^+_x \leq \|\omega\|^+_y$. Thus
\begin{align*}
\| u^{-(n+1)}(x)\|^+_x &\leq \int_x^\infty \| A(y) u^{-(n)} (y)\|^+_x \, dy \\
&\leq \int_x^\infty \|A(y) u^{-(n)} (y)\|^+_y\, dy \\
&\leq \int_x^\infty  G(y) \ \|u^{-(n)}(y)\|^+_y \, dy
\end{align*}
by \eqref{2.6}. Thus
\[
\sup_{y\geq x} \| u^{-(n+1)} (y)\|^+_y \leq \sup_{y \geq x} \| u^{-(n)}(y)\|^+_y 
\int_x^\infty G(y)\, dy
\]
proving convergence of the integral and
\[
\sup_{y\geq x} \| u^{-(n)}(y)\|^+_y \leq \bigg[ \int_x^\infty G(y)\, dy \bigg]^n
\]
inductively.

It follows that
\[
u^-(y) \equiv \sum_{n=0}^\infty u^{-(n)}(y)
\]
converges for $y\geq x_0$ where $\int_{x_0}^\infty G(y)\, dy < 1$ and that for such 
$y$, $\int_y^\infty A(w) u^- (w)\, dw$ converges and
\[
u^- (y) = \binom{1}{0} + \int_y^\infty A(w) u^- (w)\, dw
\]
so $u^-$ solves \eqref{2.1}. Since $\| u^-(y) - \binom{1}{0} \|^+_y \to 0$ as 
$y\to\infty$, we obtain (i).

Define $\tilde u^{+(n)}$ by
\begin{align*}
\tilde u^{+(0)} &= \binom{0}{1} \\
\tilde u^{+(n+1)}(x) &= \int_{x_0}^x A(y) \tilde u^{+(n)}(y)\, dy
\end{align*} 
for $x_0$ chosen so that
\begin{equation} \lb{2.7}
\int_{x_0}^\infty G(y)\, dy \leq \tfrac13.
\end{equation}
As above, using the fact that if $y<x$, then $\|\omega\|^-_x \leq \|\omega\|^-_y$ 
since $f_-$ is decreasing, we have
\[
\sup_{x\geq x_0} \|\tilde u^{+(n)}(x)\|^-_x \leq \bigg(\int_x^\infty G(y)\, dy\bigg)^n 
\leq \bigg(\frac13\bigg)^n.
\]

As in the $\|\dott\|^+$ case, we see that $\sum_{n=0}^\infty \tilde u^{+(n)} = 
\tilde u^+$ converges for $y>x_0$ and $\tilde u^+$ solves \eqref{2.1} and obeys
\[
\tilde u^+ (x) = \binom{0}{1} + \int_{x_0}^x A(y) \tilde u^+ (y)\, dy.
\]

In particular,
\[
\tilde u^+_2 (\infty) = 1 + \int_{x_0}^\infty A(y) \tilde u^+(y)\, dy 
\]
exists and $|\tilde u^+_2 (\infty) - 1| \leq \frac12$ so $\tilde u^+_2 (\infty) \equiv 
\alpha >0$. Define
\[
u^+ = \alpha^{-1} \tilde u^+
\]
and so obtain a vector-valued function $u^+_2$ with $u^+_2\to 1$ and $\abs{u^+_1 f_-}$ 
bounded. We will show that if $f_-\to 0$, then $u^+_1 f_-\to 0$. When $f_-$ does not 
go to zero, we will provide an alternative construction of $u^+$.

To prove that $u^+_1 f_- \to 0$ if $f_- \to 0$, write for $x_0 < y < x$:
\begin{equation*}
\begin{split}
u^+_1 (x) f_- (x) &= f_- (x) \int_{x_0}^y (A(w)u^+ (w))_1 \, dw \\
&\qquad + \int_y^x f_-(x) f_- (w)^{-1} f_- (w) (A(w) u^+ (w))_1 \, dw
\end{split}
\end{equation*}
so, since $f_-$ is monotone decreasing,
\begin{equation} \lb{2.8}
|u^+_2 (x) f_- (x) | \leq f_- (x) \int_{x_0}^y |(A(w) u^+(w))_1 |\, dw 
 + \int_y^\infty G(w) ||u^+ (w)\|^-_w\, dw.
\end{equation}
Given $\veps$, pick $y$ so the second integral in \eqref{2.8} is less than $\veps/2$ and 
then, since $f_- \to 0$, $x$ so that the first term is less than $\veps/2$. Thus 
$u^+_1 f_- \to 0$.

If $f_-(x)$ has a non-zero limit as $x\to\infty$, then since $f_-$ is monotone, 
$f_-(x)\geq c$. Thus, $A(x)\in L^1$, and by the same construction as used for 
$u_-$ (i.e., integrating from infinity), one can construct $u^+(x)\to\binom{1}{0}$ 
as $x\to\infty$.
\end{proof}

The situation in the discrete case is similar. \eqref{2.1} becomes
\begin{equation} \lb{2.9}
u(n+1) - u(n) = A(n) u(n).
\end{equation}
$f_\pm$ obey \eqref{2.4}, although they are only defined (and monotone) on $n=1,2,\dots$. 
$G(n)$ is defined as in \eqref{2.5}. The analog of Theorem~\ref{T2.2} holds with 
$\int_x^\infty G(y)\, dy <\infty$ replaced by 
\[
\sum_{n_0}^\infty G(n) < \infty
\]
and $x$ going through discrete values. The proof is identical with obvious changes --- 
for example, the formula for $u^{-(n+1)}$ becomes
\[
u^{-(n+1)}(j) = -\sum_{k=j}^\infty A(k) u^{-(n)}(k).
\]

\bigskip

We owe to F.~Gesztesy an illuminating remark about our result, Theorem~\ref{T2.2}, 
namely the special case $f_+f_- =1$ (which is true in some of the applications 
we will make) follows quickly from Levinson's theorem \cite{Eas,Lev}. One variant 
of Levinson's theorem says:

\begin{proposition}\lb{P2.3} Let $A$ be a $2\times 2$ matrix of the form $A_1 + A_2$ 
where $\int_{t_0}^\infty \| A_1(s)\| \, ds<\infty$, $A_2$ is diagonal with
\[
A_2 (s) = \begin{pmatrix}
\alpha_1(s) & 0 \\
0 & \alpha_2(s) \end{pmatrix}
\]
so that
\begin{equation} \lb{2.10}
\int_{t_0}^t \Real [\alpha_1(s) - \alpha_2(s)]
\end{equation}
is either bounded below or bounded above. Then, there exist solutions 
$\varphi_{1,2}$ of
\[
\dot{\varphi} = A\varphi
\]
so that
\[
\varphi_1(t)\, e^{-\int_{t_0}^t \alpha_1(s)\,ds} \to \binom{1}{0}
\]
and
\[
\varphi_2\, e^{-\int_{t_0}^t \alpha_2(s)\, ds} \to \binom{0}{1}
\]
as $t\to\infty$.
\end{proposition}

\noindent{\it Remark.} This is essentially equivalent to the general $2\times 2$ case.

\medskip

To apply this to the situation of Theorem~\ref{T2.2}, given a solution, $u$, of 
\eqref{2.2}, let $\varphi$ be defined by $\varphi_1 = u_1$, $\varphi_2 = u_2 f_+$. 
Then
\[
\varphi' = (A_1 + A_2)\varphi ,
\]
where
\[
A_1 = \begin{pmatrix}
a_{11}(x) & a_{12}(x) f_+^{-1}(x) \\
a_{21}(x) f_+(x) & a_{22}(x)
\end{pmatrix}
\]
and
\[
A_2 = \begin{pmatrix}
0 & 0 \\
0 & f'_+ / f_+(x) 
\end{pmatrix}.
\]
By hypothesis $(\int_{x_0}^\infty G(x) < \infty)$, $A_1\in L^1$ and the function 
\eqref{2.10} is $f_+(t_0)/f_+(t)$ which is bounded by $1$ (since $f_+$ is monotone). 
The two Levinson's theorem solutions obey (i) and (ii) of Theorem~\ref{T2.2}.

%%%%%%%%%%%%%%%%%%%%%%%%%%%%%%%%%%%%%%%%%%%%%%%%%%%%%%%%%%%%%
\section{Stability of Point Spectra} \lb{s3}
%%%%%%%%%%%%%%%%%%%%%%%%%%%%%%%%%%%%%%%%%%%%%%%%%%%%%%%%%%%%%

In this section, we will prove Theorems~\ref{T1.7} and \ref{T1.9}.
We will only consider the continuum case; the discrete case is similar.

\begin{proof}[Proof of Theorem~\ref{T1.7}] Since $\lambda\in P(V_0)$, $\varphi_- 
\equiv \varphi_{1, \theta(\lambda)}$ is in $L^2$ and so by hypothesis, also in $L^\infty$ 
\cite{ssg}.  Pick $\varphi_+ \equiv \varphi_{2, \theta(\lambda)}$ and use 
variation of parameters \eqref{1.21}/\eqref{1.22}. $A$ has the form \eqref{1.24}. Let 
$f_+$ be given by \eqref{1.13}. Since $\varphi_-$ is bounded,
\begin{gather*}
\abs{W(x)\varphi_+(x) \varphi_-(x)} \leq Cf_+(x) \abs{W(x)} \\
\abs{W(x)}\, \abs{\varphi_+ (x)}^2 f_+(x)^{-1} \leq f_+ (x) \abs{W(x)} \\
\abs{W(x)}\, \abs{\varphi_- (x)}^2 f_+ (x) \leq Cf_+ (x) \abs{W(x)}. 
\end{gather*}
So if $f_- (x) \equiv f_+ (x)^{-1}$, we have that $G$ given by \eqref{2.5} obeys
\[
\abs{G(x)} \leq C f_+ (x) W(x).
\]
Thus, Theorem~\ref{T2.2} is applicable, so there is a solution, $\psi$, of the 
perturbed Schr\"odinger equation of the form:
\begin{equation} \lb{3.1}
\psi(x) = u^-_1 (x) \varphi_-(x) + u^-_2 (x) \varphi_+(x)
\end{equation}
with $u^-_1$ bounded and with $u^-_2 (x) f_+(x)$ bounded. Since $\varphi_- 
\in L^2$, $u^-_1 \varphi_- \in L^2$. Moreover, since $u^-_2 f_+$ is bounded, 
\eqref{1.13} says that
\[
\abs{u^-_2 (x) \varphi_+(x)} \leq C(1 + \abs{x})^{-1}
\]
which is also in $L^2$. Thus $\psi\in L^2$. 
\end{proof}

Corollary~\ref{C1.8} follows immediately since $f_+$ is exponentially bounded. 
Results of this genre are well known; see, for example, Hartman \cite{Har}. 
We proceed to prove Theorem~\ref{T1.9}:

\begin{proof}[Proof of Theorem~\ref{T1.9}] If $\gamma=0$, any solution, $\varphi$, 
of the unperturbed equation has
\begin{equation} \lb{3.2}
\abs{\varphi(x)} + \abs{\varphi'(x)} \leq C_\veps \, e^{\veps\abs{x}}
\end{equation}
so, by hypothesis, the $A$ of \eqref{1.24} is in $L^1$ for any choice of 
$\varphi_\pm$. Thus, by standard theory (or Theorem~\ref{T2.2} with $f_+ = 
f_- = 1$), any solution $\psi$ also obeys \eqref{3.2} which implies that 
$\gamma = 0$.

Now suppose that $\gamma >0$. By the Ruelle-Osceledec Theorem \cite{Ru}, there is a 
solution $\varphi_-(x)$ ($\equiv \varphi_{1, \theta(\lambda)}$) for the $V_0$ 
equation with
\[
\lim_{x \rightarrow \infty}\left[
\f{1}{\abs{x}} \, \ln [\abs{\varphi_-(x)}^2 + \abs{\varphi'_- (x)}^2] \right] = -\gamma.
\]
Any linearly independent solution and, in particular, $\varphi_+ = 
\varphi_{2, \theta(\lambda)}$ obeys
\[
\lim_{x \rightarrow \infty} \left[
\f{1}{\abs{x}} \, \ln [\abs{\varphi_+ (x)}^2 + \abs{\varphi'_+(x)}^2]\right] = \gamma.
\]
In particular, for any $\veps_1 >0$,
\begin{equation} \lb{3.3}
\abs{\varphi_+(x)} \leq C_{\veps_1} e^{(\gamma + \veps_1)\abs{x}}, 
\qquad \abs{\varphi_-(x)} \leq C_{\veps_1} \, e^{-(\gamma-\veps_1)\abs{x}}.
\end{equation}
Pick $f_\pm (x) = e^{(\pm 2\gamma + 2\veps_1)\abs{x}}$ where $\veps_1$ is chosen so 
that $\veps_1 <\gamma$ (so $f_-$ is decreasing) and $\veps_1 < \f14 \veps$ where 
$\veps$ is given in the hypothesis of the theorem. By the estimates of \eqref{3.3}, 
\[
\abs{G(x)} \leq e^{4\veps_1 \abs{x}} \abs{W(x)}
\]
so $G\in L^1$ since $4\veps_1 <\veps$. Theorem~\ref{T2.2} applies and we get 
solutions $\psi_\pm$ of the perturbed equation with
\begin{align*}
\abs{\psi_+ - u^+_2 \varphi_+} &\leq u^+_1 \abs{\varphi_-} \\
&\leq C_{\veps_1} u^+_1 f_- \, e^{(\gamma - \veps_1)\abs{x}}
\end{align*}
and a similar estimate for $\psi'_+$. It follows that
\[
\f{\| \psi_+\|_x}{\| \varphi_+\|_x} - 1 \to 0
\]
as $x\to \infty$ where $\| g\|_x = (\abs{g(x)}^2 + \abs{g'(x)}^2)^{1/2}$. 
Similarly $\|\psi_- \|_x / \| \varphi_- \|_x \to 1$. Thus not only is the Lyapunov 
exponent the same, but even the subexponential corrections are unchanged. 
\end{proof}

%%%%%%%%%%%%%%%%%%%%%%%%%%%%%%%%%%%%%%
\section{Power Law Theorems} \lb{s4}
%%%%%%%%%%%%%%%%%%%%%%%%%%%%%%%%%%%%%%%

In this section we will prove the following result that essentially includes 
Theorems~\ref{T1.10} and \ref{T1.11} as corollaries. (We will need to make 
an additional argument for $\beta = 0$.)

\begin{theorem} \lb{T4.1} In the continuum case, suppose $V_0$ is bounded. Let 
$\lambda\in S_0 (V_0)$ be a regular energy with $\beta (\lambda, V_0) >0$. Suppose 
that
\[
\abs{W(x)} \leq C (1 + \abs{x})^{-\eta}
\]
for some $\eta > \beta(\lambda, V_0)^{-1}$. Then $\lambda\in S (V_0 + W)$ and 
\[
\beta (\lambda, V_0 + W) = \beta (\lambda, V_0).
\]
\end{theorem}

Interestingly enough, we will apply Theorem~\ref{T2.2} in a situation where $f_+ f_- 
\neq 1$ but is strictly bigger. Essentially, we will not want to take $f_-$ as small as 
$f^{-1}_+$ because we will need the error estimate $u^+_1 f_-\to 0$ to be stronger 
than $u^+_1 f^{-1}_+ \to 0$.

To employ the ideas of Jitomirskaya-Last, we need to relate estimates involving an 
integral of a product of $\varphi_+, \varphi_-, W$ and $1, f_+$ or $f_-$ to 
$\| \dott\|_L$. The following is useful:

\begin{lemma} \lb{L4.2} If
\begin{equation} \lb{4.1}
\abs{Q(x)} \leq C_1 (1 + \abs{x})^{-a}
\end{equation}
and
\begin{equation} \lb{4.2}
\| \varphi_+ \|_L \| \varphi_-\|_L \leq C_2 (1 + L)^b
\end{equation}
and $a>b$, then
\[
\int_0^\infty \abs{Q(x) \varphi_+ (x) \varphi_- (x)} \, dx < \infty.
\]
\end{lemma}

\begin{proof} Let $g(x) = \int_0^x \abs{\varphi_+ (y) \varphi_-(y)}\, dy$. 
By the Schwarz inequality and \eqref{4.2},
\begin{equation} \lb{4.3}
\abs{g(x)} \leq C_2 (1 + \abs{x})^b
\end{equation}
and, of course,
\begin{equation} \lb{4.4}
g(0) = 0.
\end{equation}
Use \eqref{4.1} to write
\begin{align*}
\int_0^L  |Q(x) \varphi_+ (x) & \varphi_-(x)|\, dx \\
&\leq C_1 \int_0^L (1 + \abs{x})^{-a} \, \f{dg}{dx}\, dx \\
&= C_1 a \int_0^L (1 + \abs{x})^{-a-1} g(x) \, dx + C_1 ( 1 + \abs{L})^{-a} g(L).
\end{align*}
There is no boundary term at $x=0$ by \eqref{4.4}.

Now use \eqref{4.3} and $b < a$ to see
\begin{align*}
\lim_{L\to\infty} \int_0^L \abs{Q(x) \varphi_+ (x) \varphi_-( x)}\,dx & \leq 
C_1 C_2 \, a \int_0^\infty (1 + \abs{x})^{b-a-1}\, dx \\
&< \infty.
\end{align*}
\end{proof}

The next step is obtaining power-law upper and lower bounds on $\| \varphi_\pm \|_L$. 
In principle, the upper and lower powers could be different with oscillation between 
the two powers of growth.

\begin{lemma} \lb{L4.3} Let $\lambda$ be a regular energy with $\lambda\in S$ and 
$\beta (\lambda) > 0$. Let $\varphi_- = \varphi_{1, \theta(\lambda)}$ and 
$\varphi_+ = \varphi_{2, \theta(\lambda)}$. Then for any $\veps > 0$, there are 
constants $C_1, C_2, C_3, C_4$ {\rm(}$\veps$-dependent{\rm)} so that for $L$ large,
\begin{align} 
C_2 L^{1- 1/2 \beta - \veps} & \leq \|\varphi_- \|_L  \leq C_1 L^{1/2 + \veps} \lb{4.5} \\
C_4 L^{1/2 - \veps} & \leq \| \varphi_+ \|_L \leq C_3 L^{1/2 \beta + \veps}. \lb{4.6}
\end{align}
\end{lemma}

\begin{proof} The definition of regularity says \eqref{1.15} which is the $C_1$ 
estimate in \eqref{4.5}. \eqref{1.16} then implies the $C_4$ estimate in \eqref{4.6}.

By \eqref{1.11}, if $\tilde\beta < \beta$, then $\|\varphi_-\|_L \geq 
\|\varphi_+\|^{\tilde\beta}_L$ for $L$ large which, given the $C_1$ estimate, implies 
the $C_3$ estimate in \eqref{4.6}. Using \eqref{1.16} again, we get the $C_2$ 
estimate in \eqref{4.5}. 
\end{proof}

At first sight, it might appear that all one needs on $\| \varphi_\pm \|_L$ are upper 
bounds because they are all that enter in proving the applicability of Theorem~\ref{T2.2}. 
But one wants to apply Theorem~\ref{T2.2} to show that
\begin{equation} \lb{4.7}
\f{\|\psi_+\|_L}{\|\varphi_+\|_L}\to 1, \qquad \f{\|\psi_-\|_L}{\|\varphi_-\|_L} 
\to 1
\end{equation}
as $L\to\infty$. Consider the second part of \eqref{4.7}. We have
\[
\psi_- = u^-_1 \varphi_- + u^-_2 \varphi_+.
\]
Since $u^-_1 \to 1$, we have that
\[
\f{\abs{ \|\psi_- \|_L - \|\varphi_- \|_L}}{\| \varphi_- \|_L} \leq 
\f{\|\psi_- -\varphi_- \|_L}{\|\varphi_- \|_L} \leq 
\f{\| \psi_-  -u^-_1 \varphi_- \|_L}{\|\varphi_- \|_L} + o(1)
\]
and so it is natural to prove the desired relation by showing
\[
\f{\|u^-_2 \varphi_+\|_L}{\|\varphi_-\|_L} \to 0.
\]
All we basically know about $u^-_2$ is $f_+ u^-_2 \to 0$. Thus

\begin{lemma}\lb{L4.4} Suppose $G(x)\in L^1$. In order for \eqref{4.7} to hold, it 
suffices that for large $L$, 
\begin{equation} \lb{4.8}
\f{\|\varphi_+\|_L}{\| \varphi_- \|_L} \leq Cf_+(L)
\end{equation}
and
\begin{equation} \lb{4.9}
\f{\|\varphi_-\|_L}{\| \varphi_+\|_L} \leq Cf_- (L).
\end{equation}
\end{lemma}

By \eqref{4.5}/\eqref{4.6}, we have \eqref{4.8} if $f_+(L)=L^{\mu_+}$ with
\[
\f{1}{2\beta} - \bigg( 1 - \f{1}{2\beta}\bigg) < \mu_+.
\]

If we only apply the similar bound for \eqref{4.9}, we see that we need $f_- (L) \geq 
L^{2\veps}$ which is incompatible with $f_-$ decreasing. We therefore do not gain 
from \eqref{4.5}/\eqref{4.6} and instead define $f_-(x) \equiv 1$ so that \eqref{4.9} 
holds since $\| \varphi_-\|_L$ is subordinate. Thus we will take
\begin{equation} \lb{4.10}
f_+(x) = x^{\mu_+}, \qquad f_-(x) = 1
\end{equation}
with
\begin{equation} \lb{4.11}
\mu_+ > \f{1}{\beta} - 1.
\end{equation}

\begin{proof}[Proof of Theorem~\ref{T4.1}] By the above analysis, if we take $f_+, f_-$ 
to obey \eqref{4.10}/\eqref{4.11}, we have \eqref{4.7} so long as Theorem~\ref{T2.2} is 
applicable. But \eqref{4.7} implies that $[ \|\psi_-\|_L/ \|\psi_+\|^{\tilde\beta}_L] 
\big/ [ \|\varphi_- \|_L / \| \varphi_+\|^{\tilde\beta}_L] \to 1$ and thus by
Proposition~\ref{P1.3}, $\beta(\lambda, V_0+W) = \beta (\lambda, V_0)$.

To apply Theorem~\ref{T2.2}, we need $G$ to be in $L^1$. By Lemma~\ref{L4.2} and the 
upper bounds in \eqref{4.5}/\eqref{4.6}, this is true if the following three 
inequalities hold
\begin{alignat}{2}
\f12 + \f{1}{2\beta} &< \eta \quad && \longleftarrow (a_{11}, a_{22}\text{ terms}) \lb{4.12} \\
1 + \, \mu_+ &< \eta \quad &&\longleftarrow (a_{21}\text{ terms}) \lb{4.13} \\
\f{1}{\beta} &< \eta \quad && \longleftarrow (a_{12}\text{ terms}). \lb{4.14}
\end{alignat}

By the basic hypothesis of the theorem, $\eta >\beta^{-1}$ and, of course, $\beta^{-1} 
\geq 1$. Thus \eqref{4.12} and \eqref{4.14} hold, and to get \eqref{4.13} and \eqref{4.11}, 
we need only choose $\mu_+ >0$ so that
\[
\f{1}{\beta} < 1 + \mu_+ < \eta.
\]
This can be done since $\beta^{-1} \geq 1$.
\end{proof}

Theorem~\ref{T1.10} is an immediate corollary of Theorem~\ref{T4.1} as is 
Theorem~\ref{T1.11} if $\beta(\lambda, V_0) \neq 0$. In case $\beta(\lambda, V_0)=0$, 
then we claim $\beta (\lambda, V_0 + W) =0$ for if not, we can turn this argument 
around (think of $V_0 = (V_0 + W) - W$) and find that $\beta (\lambda, V_0) = \beta
(\lambda, V_0 +W)\neq 0$. That means $\beta(\lambda, V_0 + W)=0$ which implies 
$\lambda\in P(V_0 + W)$ or $\lambda \in S(V_0 + W)$ with $\beta =0$.

The condition $\eta > \beta^{-1}$ of Theorem~\ref{T4.1} is needed because we assume 
no extra information about the behavior of $\| u_1\|_L$ and $\|u_2\|_L$ other than 
the value of $\beta$. If one has additional information, one can often do better. 
Here is an extreme example, but one that holds in some explicit examples.

\smallskip
\noindent{\bf Definition.} We say there is power Lyapunov-Osceledec
 behavior with 
exponent $\gamma >0$ at energy $\lambda$ if and only if there exist solutions 
$\varphi_{1,\theta(\lambda)}$ and $\varphi_{2,\theta(\lambda)}$ with 
\begin{align*}
\lim_{x\to\infty} \f{\ln [ \abs{\varphi_{1,\theta(\lambda)}(x)}^2 + 
\abs{\varphi'_{1,\theta(\lambda)}(x)}^2]}{\ln \abs{x}} &= - \gamma \\
\lim_{x\to\infty} \f{\ln [ \abs{\varphi_{2,\theta(\lambda)}(x)}^2 + 
\abs{\varphi'_{2,\theta(\lambda)}(x)}^2]}{\ln \abs{x}} &=  \gamma.
\end{align*}

\smallskip
\noindent{\it Notes.} 1. In the discrete case, replace $\varphi'(x)$ by 
$\varphi (n+1)$.

\smallskip
2. Under these circumstances, if $\gamma < 1/2$ and $V$ is bounded, we have 
$\| \varphi_1\|_L \sim L^{-\gamma + 1/2}$,
 $\|\varphi_2\|_L \sim L^{\gamma + 1/2}$ 
(where $\sim$ means up to factors of $L^\veps$) so $\lambda\in S(V_0)$ and 
$\beta(\lambda) = (1/2 - \gamma)/(1/2 + \gamma)$ and $\alpha(\lambda) = 
1-2\gamma$.

\smallskip
3. One example where it is known \cite{KLS} there is power Lyapunov-Osceledec 
behavior is the discrete $n^{-1/2}$ decaying Anderson
 model where $V_\omega (n) 
= \lambda n^{-1/2} X_\omega (n)$ where the $X_\omega$
 are bounded i.i.d.'s with 
$E(X_\omega)=0$, $E(X^2_\omega)=1$. Then there is power Lyapunov-Osceledec 
behavior with $\gamma = \lambda^2/(8-2E^2)$ in the region $\abs{E} < 2$.

\smallskip
\begin{theorem} \lb{T4.5} Suppose $V_0$ has power
 Lyapunov-Osceledec behavior with 
$\gamma > 0$ at energy $\lambda$ and that 
\[
\int (1 + \abs{x})^\veps \abs{W(x)}\, dx < \infty
\]
for some $\veps > 0$. Then $V_0 + W$ has power
 Lyapunov-Osceledec behavior at energy 
$\lambda$ with the same value of $\gamma$.
\end{theorem}

The proof is essentially identical to the proof of
 Theorem~\ref{T1.9} (in Section~\ref{s3})  
with $f_\pm (x) = (1+ \abs{x})^{2(\pm\gamma + \veps)}$.

This shows the improvement over the power in Theorem~\ref{T4.1}.
 Instead of $\eta > 
\beta^{-1}$, we only need $\eta > 1$.

%%%%%%%%%%%%%%%%%%%%%%%%%%%%%%%%%%%%%%%
\section{Appendix: WKB Asymptotic Behavior} \lb{a1}
%%%%%%%%%%%%%%%%%%%%%%%%%%%%%%%%%%%%%%%

In this appendix, we illustrate with an example how
Theorem~\ref{T2.2} can be used to obtain precise 
asymptotic behavior of solutions in some concrete situations 
(where the perturbation does not even have to be decaying).
Namely, we  show how to use Theorem~\ref{T2.2} to prove the existence of 
WKB solutions at $+\infty$ for
\begin{equation} \lb{3.4}
-\psi'' + V\psi = \lambda\psi,
\end{equation}
when $V=V_1 + V_2$ with
\begin{equation} \lb{3.5}
V_1 \in L^1, \quad V'_2 \in L^1, \quad V_2 (x) \to 0 \qquad\text{as } x\to\infty,
\end{equation}
when $\lambda\neq 0$. For $\lambda >0$, it is well known that such solutions exist 
(see, e.g., \cite{spams}).  
For $\lambda <0$, 
one can also apply Levinson's theorem (see \cite{Lev} or \cite{CL}, 
Theorem 8.1) to prove this result, but it is nice to get it from Theorem~\ref{T2.2}.
 As a preliminary, we note one can try 
an Ansatz, \eqref{1.21}/\eqref{1.20} for solving \eqref{3.4} even if $\varphi_\pm$  
do not solve a related Schr\"odinger equation. The result is that $u$ still obeys 
\eqref{1.23} but $A$ is now given by
\begin{equation} \lb{3.6}
A(x) = w(x)^{-1}\begin{pmatrix}
-\varphi_+(x) (H_\lambda \varphi_-)(x) & -\varphi_+(x)(H_\lambda \varphi_+)(x) \\
\varphi_-(x) (H_\lambda \varphi_-)(x) & \varphi_-(x) (H_\lambda \varphi_+)(x) 
\end{pmatrix},
\end{equation}
where
\begin{equation} \lb{3.7}
w(x) = \varphi_-(x) \varphi'_+(x) - \varphi_+(x) \varphi'_-(x)
\end{equation}
and $H_\lambda$ is the differential expression
\begin{equation} \lb{3.8}
H_\lambda = -\f{d^2}{dx^2} + V-\lambda.
\end{equation}

We can now prove

\begin{theorem} \lb{T3.1} Let $V$ obey \eqref{3.5} and $\lambda \neq 0$. If 
$\lambda < 0$, let
\[
\varphi_\pm (x) = \exp (\pm \eta(x)),
\]
where
\[
\eta(x) = \int_{s_0}^x \sqrt{-\lambda + V_2(s)}\, ds,
\]
and $s_0$ is chosen so that $\abs{V_2 (s)} \leq \abs{\lambda}$ for $s > s_0$.
If $\lambda >0$, let
\[
\varphi_\pm = \exp (\pm\eta (x)),
\]
where 
\[
\eta = i \int_{s_0}^x \sqrt{\lambda - V_2(s)}\, ds.
\]

Then there exist solutions $\psi_\pm$ of \eqref{3.4} so
\begin{align*}
\psi_\pm(x) &= \varphi_\pm (x) \, (1 + o(1)) \\
\psi'_\pm (x) &= \varphi'_\pm (x) \, (1 + o(1)) 
\end{align*}
as $x\to\infty$.
\end{theorem}

\begin{proof} Consider first the case $\lambda <0$. Then $\varphi'_\pm = \pm \eta' 
e^{\pm\eta}$ and $\varphi''_\pm = (\pm \eta'' + (\eta')^2)\varphi_\pm$ and thus,  
since $(\eta')^2 = -\lambda + V_2$,
\begin{align*}
H_\lambda \varphi_\pm &= [\pm\eta'' + V_1] \varphi_\pm \\
&= \bigg[ \f{\pm V'_2}{2\eta'} + V_1 \bigg] \varphi_\pm
\end{align*}
so we define
\[
Q_\pm = \f{\pm V'_2}{2\eta'} + V_1.
\]
Since $\eta' \to \sqrt{-\lambda}$ as $x\to \infty$, we see that $Q_\pm \in L^1$. 
Moreover,
\begin{equation} \lb{3.9}
w(x) = 2\eta' \to 2\sqrt{-\lambda} \qquad\text{as } x\to\infty.
\end{equation}
It follows with $f_\pm = \varphi^2_\pm$ (so $f_+ = f^{-1}_-)$ and $A$ given by 
\eqref{3.6} that $G(x)\in L^1$ since $Q_\pm(x)\in L^1$.

Applying Theorem~\ref{T2.2}, there are solutions $\psi_\pm$ with
\begin{align*}
\psi_\pm &= \varphi_\pm (1 + o(1)) + \varphi_\mp (f^{-1}_\mp) o(1) \\
&= \varphi_\pm (1 + o(1)) 
\end{align*}
and similarly for $\psi'_\pm$.

The calculation for $\lambda >0$ is similar, except we use $\abs{\varphi_+} = 
\abs{\varphi_-} =1$ in that case to pick $f_+ = f_- =1$.
\end{proof}

In this paper, we considered only perturbations which are absolutely integrable.
It is reasonable to ask what one can expect for stronger perturbations, 
for example, in situations where there is Lyapunov behavior. 
While in general the picture is not complete, we provide 
a sample result which gives $L^2$ stability under additional 
assumptions on the behavior of solutions of the unperturbed equation.
As a bonus, we also obtain a stronger version of Theorem~\ref{T3.1}
in the case $\lambda<0.$ 
%In particular, the Hartman-Wintner 
%theorem (see, e.g. \cite{Eas})
%shows that a result similar to
% Theorem~\ref{T3.1} on the WKB behavior of solutions 
%remains true for $\lambda<0$ and $V_1 \in L^2, V_2=0.$ 
%Even in more generality, we have 
\begin{proposition}\label{l2}
 Assume that there exist functions $\varphi_\pm(\lambda,x)$
such that $(-\frac{d^2}{dx^2}+V-\lambda) \varphi_\pm = U_\pm \varphi_\pm,$
with $U_\pm \in L^2,$ and that the inverse of the Wronskian 
$W[\varphi_-, \varphi_+]^{-1}$
is bounded. Define functions
\begin{equation}\label{corfun}
\eta_\pm(\lambda,x) = \varphi_\pm e^{\pm
\int_0^x  \frac{U_\pm \varphi_+ \varphi_-}{W[\varphi_-,\varphi_+]} \,dt}
\end{equation}
and the kernel 
\begin{equation}\label{ker}
K(x,y) = \varphi^2_+(x) \varphi^2_-(y) e^{-\int_x^y
\frac{(U_-+U_+)\varphi_-\varphi_+}{W[\varphi_-,\varphi_+]} \,dt}.
\end{equation}
Assume in addition that 
\begin{equation}\label{bigger}
{\rm inf}_x |\varphi_+(x) \varphi_-(x)| 
\geq c>0
\end{equation}
and
\begin{equation} \label{kercon}
\int\limits_0^\infty {\rm sup}_x |K(x,x+y)| \,dy 
< \infty, \,\,\,{\rm sup}_{y \geq x}|K(x,y)| \leq C.
\end{equation}
Then there exist solutions $\psi_\pm$ of the equation 
$(H+V-\lambda)\psi_\pm =0$ with the asymptotic behavior 
\begin{equation}\label{asbe}
\psi_\pm(\lambda,x) =  \eta_\pm(\lambda,x) (1+o(1)).  
\end{equation}
\end{proposition}
\noindent {\it Remarks.} 1. In order for \eqref{kercon} to hold, 
one needs, roughly speaking, Lyapunov behavior at $\lambda$ and 
moreover $\phi_+(x,\lambda) \phi_-(x,\lambda) \sim \text{const}$ 
(or grows very slowly) for large $x.$ 

\smallskip
2. In the case where $V_2=0,$ the result follows from 
 the Hartman-Wintner theorem (see, e.g., \cite{Eas}).
%shows that a result similar to
% Theorem~\ref{T3.1} on the WKB behavior of solutions 
%remains true for $\lambda<0$ and $V_1 \in L^2, V_2=0.$ 

\smallskip
3. Notice that the asymptotic behavior of solutions of the perturbed
equation differs from $\varphi_\pm$ by an additional factor. 

\medskip
Before sketching the proof, let us illustrate the result 
with the following generalization of Theorem~\ref{T3.1}
for $\lambda<0$.
\begin{corollary}
%The result of Theorem~\ref{3.1} remains true under the 
Assume that $V=V_1+V_2,$ $V_1 \in L^2,$ $V_2' \in L^2,$ 
$V_2(x) \rightarrow 0$ as $x \rightarrow \infty.$
Then for $\lambda <0$ there exist solutions $\psi_\pm$ of the equation 
$-\psi'' +V \psi = \lambda \psi$ such that
\[ \psi_\pm (x,\lambda) = \eta_\pm (x,\lambda)
(1+o(1)). \]
Here $\eta(x,\lambda)$ is given by \eqref{corfun} with 
\[ \varphi_\pm(x,\lambda) = \exp \left(
\int_{s_0}^x \sqrt{-\lambda +V_2(s)}\,ds \right), \]
and $s_0$ is such that $|V_2(s)| < |\lambda|$ for $s>s_0.$    
\end{corollary} 
\noindent \it Remark. \rm  For $\lambda>0,$ the result is generally not true.
It holds for a.e. $\lambda>0$ for $V_1 \in L_p,$ $V'_2
\in L^p$ with $p<2$ \cite{CK1}. It is not known if 
the result remains true for $p=2$ and a.e. $\lambda>0.$  
\begin{proof}
Choosing $\varphi_\pm$ as in the statement of the corollary, 
one directly 
verifies that all conditions of Proposition~\ref{l2} hold. 
Notice that the Lyapunov behavior is preserved, since 
\begin{gather*}
 W[\phi_-,\phi_+] = -2\sqrt{-\lambda +V_2}, \\ 
 U_\pm = V_1 \mp \frac{V'_2}{2\sqrt{-\lambda +V_2}} 
\end{gather*}
and therefore the additional factor in \eqref{corfun}
is bounded by $e^{Cx^{1/2}}$.  
\end{proof}

We now sketch the proof of Proposition~\ref{l2}. 
Seeking solution $\psi(x)$ of the equation $-\psi''+V\psi = \lambda \psi, $
apply variation of parameters-type transformation 
\[ \left( \begin{array}{c} \psi \\ \psi' \end{array} \right) =
 \left( \begin{array}{cc} 
\varphi_- & \varphi_+ \\ \varphi'_- &  \varphi'_+ \end{array} \right) u(x), \]
obtaining a system 
\[ u'(x) = \frac{1}{W[\varphi_-, \varphi_+]} \left( \begin{array}{cc} 
-U_- \varphi_- \varphi_+ & -U_+ \varphi_+^2 \\ U_- \varphi_-^2
 & U_+ \varphi_-\varphi_+ 
\end{array} \right) u(x). \]
Do one more transformation to bring this system to a simpler form:
\[ u(x) = \left( \begin{array}{cc} e^{- \int_0^x 
\frac{U_- \varphi_-\varphi_+}{W[\varphi_-, \varphi_+]} \,dt}
& 0 \\ 0 &  e^{ \int_0^x \frac{
U_+ \varphi_-\varphi_+}{W[\varphi_-, \varphi_+]} \,dt} 
\end{array}  \right) z(x), \]
then 
\begin{equation}\label{finsys}
z'(x) =  \left( \begin{array}{cc} 
0 & -\frac{U_+ \varphi_+^2}{W[\varphi_-, \varphi_+]} e^{ \int_0^x 
\frac{(U_-+U_+) \varphi_-\varphi_+}{W[\varphi_-, \varphi_+]} \,dt} \\
-\frac{U_- \varphi_-^2}{W[\varphi_-, \varphi_+]} e^{-
 \int_0^x \frac{(U_-+U_+) \varphi_-\varphi_+}{W[\varphi_-, \varphi_+]}
 \,dt} & 0 
\end{array} \right)z(x).
\end{equation}
One can obtain the formal series for solutions of \eqref{finsys} by 
iteration; starting 
with the vector $(1,0)^T$ will lead to the solution $\psi_-(x).$ 
Properties \eqref{bigger} and \eqref{kercon} allow one to prove
the convergence of this series and \eqref{asbe} using elementary estimates. 
We leave the details to the interested reader.  

\bigskip
\noindent{\bf Acknowledgments.}  A.~Kiselev and Y.~Last would like to thank
 T.~Tombrello for 
the hospitality of Caltech. We would like to thank F.~Gesztesy and E.~Harrell 
for useful comments.

\bigskip

%%%%%%%%%%%%%%%%%%%%%%%%%%%

\end{document}